\documentclass[12pt]{amsart}
\usepackage{amsfonts}

\theoremstyle{definition}

\theoremstyle{remark}

\newcommand{\const}{\mathop{\rm const}\limits}

\newcommand{\supp}{\mathop{\rm supp}\limits}

\begin{document}

\begin{center}

{\bf WEIGHTED HARDY-LITTLEWOOD AVERAGE OPERATORS\\
 ON BILATERAL GRAND LEBESGUE SPACES} \par

\vspace{3mm}
{\bf E. Ostrovsky}\\

e-mail: galo@list.ru \\

\vspace{3mm}

{\bf L. Sirota}\\

e-mail: sirota3@bezeqint.net \\

\vspace{3mm}

Department of Mathematics and Statistics, Bar-Ilan University, 59200, Ramat Gan, Israel. \\

\vspace{3mm}

\begin{center}
 Abstract. \\

\end{center}

{\it We obtain in this short article the non-asymptotic exact estimations for the norm of (generalized)
weighted Hardy-Littlewood average integral operator in the so-called
Bilateral Grand Lebesgue Spaces. We also give examples to show the
sharpness of these inequalities.} \\

\end{center}

\vspace{3mm}

2000 {\it Mathematics Subject Classification.} Primary 37B30,
33K55; Secondary 34A34, 65M20, 42B25.\\

\vspace{3mm}

Key words and phrases: norm, Grand and ordinary Lebesgue-Riesz Spaces, weighted 
Hardy-Littlewood average operator,
exact estimations.\\

\vspace{3mm}

\section{Introduction. Statement of problem.}

\vspace{3mm}

"In 1984, Carton-Lebrun, Fosset \cite{Fosset1} defined the so-called weighted Hardy-Littlewood average operator
$  U_{\phi}. $ Let $ \phi: [0, 1] \to [0,\infty) $ be a non-negative non-zero measurable function. Let $ f $ be a measurable
real (or complex) valued function on $ R^d, \ d = 1,2, \ldots; $ one then defines the weighted Hardy–Littlewood average $U_{\phi}[f](\cdot) $ as
follows:

$$
U_{\phi}[f](x) = \int_0^1 f(tx) \ \phi(t) \ dt, \ x \in R^d. \eqno(1.1)
$$

 In this view, we think it may be suitable to call the operator $  U_{\phi} $ the {\it generalized Hardy
operator.} \par

 In \cite{Xiao1} J.Xiao obtained that the generalized Hardy operator $ U_{\phi} $ is bounded on
 $ L_p(R^d),  1 \le p \le \infty, $ if and only if

 $$
 \theta(p) = \theta_d(p)  \stackrel{def}{=} \int_0^1 t^{-d/p} \ \phi(t) \ dt < \infty. \eqno(1.2)
 $$

Moreover,

$$
||U_{\phi}||\{L(p) \to L(p) \} = \theta(p). \eqno(1.3)
$$
 It was investigated also the boundedness of this operator in the BMO spaces etc." \cite{Zhao1}.\par

 Here as usually

 $$
|f|L(p) = |f|_p = \left[ \int_{R^d} |f(x)|^p \ dx \right]^{1/p}.
 $$

 The results seems to be of interest as it is related closely to the Hardy integral inequality.
For example, if $ \phi = 1 $ and $ d = 1, $ then $ U_{\phi} $ coincides to the classical Hardy operator.

\vspace{3mm}

{\bf  Our aim is to extend the results of J.Xiao about  boundedness of Hardy-Littlewood operator on
the so-called  (bilateral)  Grand Lebesgue Spaces (GLS).}\par
 We intend to obtain the exact value of  correspondent multiplicative constants.\par

\vspace{3mm}

 The norm estimates for different operators acting in GLS  see, e.g. in the works
\cite{Liflyand1}, \cite{Ivaniec2}, \cite{Ostrovsky2}, \cite{Ostrovsky100}, \cite{Ostrovsky101},
\cite{Ostrovsky3}, \cite{Ostrovsky4}, \cite{Ostrovsky6}, \cite{Ostrovsky7}. \par

 Some application in the theory of non-linear PDE  are described  in
\cite{Ostrovsky5}, \cite{Ostrovsky8}.\par

\vspace{3mm}

 We use symbols $C(X,Y),$ $C(p,q;\psi),$ etc., to denote positive
constants along with parameters they depend on, or at least
dependence on which is essential in our study. To distinguish
between two different constants depending on the same parameters
we will additionally enumerate them, like $C_1(X,Y)$ and
$C_2(X,Y).$ The relation $ g(\cdot) \asymp h(\cdot), \ p \in (A,B), $
where $ g = g(p), \ h = h(p), \ g,h: (A,B) \to R_+, $
denotes as usually

$$
0< \inf_{p\in (A,B)} h(p)/g(p) \le \sup_{p \in(A,B)}h(p)/g(p)<\infty.
$$
The symbol $ \sim $ will denote usual equivalence in the limit
sense.\par
We will denote as ordinary the indicator function
$$
I(x \in A) = 1, x \in A, \ I(x \in A) = 0, x \notin A;
$$
here $ A $ is a measurable set.\par
 All the passing to the limit in this article may be grounded by means
 of Lebesgue dominated convergence theorem.\par

\bigskip

\section{Grand Lebesgue Spaces (GLS). }

\vspace{3mm}

We recall first of all here  for reader conventions some definitions and facts from
the theory of GLS spaces.\par

Recently, see
\cite{Fiorenza1}, \cite{Fiorenza2},\cite{Ivaniec1}, \cite{Ivaniec2}, \cite{Jawerth1},
\cite{Karadzov1}, \cite{Kozatchenko1}, \cite{Liflyand1}, \cite{Ostrovsky1}, \cite{Ostrovsky2} etc.
 appear the so-called Grand Lebesgue Spaces GLS
 $$
 G(\psi) = G = G(\psi ; A;B);  \ A;B = \const; \ A \ge 1, \ B \le \infty
 $$
spaces consisting on all the measurable functions $ f : X \to R  $ with finite norms

$$
||f||G(\psi) \stackrel{def}{=} \sup_{p \in (A;B)} \left[\frac{|f|_p}{\psi(p)} \right]. \eqno(2.1)
$$

 Here $ \psi = \psi(p), \ p \in (A,B) $ is some continuous positive on the {\it open} interval $ (A;B) $ function such
that

$$
\inf_{p \in(A;B)} \psi(p) > 0. \eqno(2.2)
$$

We will denote
$$
\supp(\psi) \stackrel{def}{=} (A;B).
$$

The set of all such a functions with support $ \supp(\psi) = (A;B) $ will be denoted by  $  \Psi(A;B). $  \par

This spaces are rearrangement invariant; and are used, for example, in
the theory of Probability, theory of Partial Differential Equations,
 Functional Analysis, theory of Fourier series,
 Martingales, Mathematical Statistics, theory of Approximation  etc. \par

 Notice that the classical Lebesgue-Riesz spaces $ L_p $  are extremal particular case of Grand Lebesgue Spaces, see
 \cite{Ostrovsky2},  \cite{Ostrovsky100}. \par

\bigskip

\section{Main result: upper and lower  estimations for Hardy - Littlewood average operator}

\vspace{3mm}

 Let the function $ \phi = \phi(t)  $ with it the function $ \theta = \theta(p) $ be a given.
Let also $ \psi = \psi(p) $ be any function described below.  Assume that

$$
\supp \psi \cap \supp \theta  =: (A,B)   \ne \emptyset,  \eqno(3.1)
$$
and denote

$$
\psi_{\theta}(p) = \psi(p) \cdot \theta(p), \ p \in (A,B).
$$

\vspace{3mm}

{\bf Theorem 3.1. }

$$
||U_{\phi} [f]||G\psi_{\theta} \le 1 \cdot ||f||G\psi, \eqno(3.2)
$$
\vspace{3mm}
{\it where the constant "1" is the best possible. } \par

\vspace{3mm}

{\bf Proof.} The {\it upper estimate}  may be proved very simple. Indeed, let $ f \in G \psi, \ f \ne 0. $ Then we have
by definition of the $ G \psi  $ norm $ |f|_p \le ||f||G\psi \cdot \psi(p). $ We deduce using Xiao inequality:

$$
|U_{\phi} [f]|_p \le |f|_p \cdot \theta(p) \le ||f||G\psi \cdot \psi(p) \theta(p) = ||f||G\psi \cdot \psi_{\theta}(p), \eqno(3.3)
$$
which implies (3.2). \par
\vspace{3mm}

{\it Lover bound.} Let the value $  \epsilon $ and the function $ \phi $ be a fix.  Denote $ \Delta = 1/\epsilon, $

$$
\theta_{\epsilon}(p) = \int_{\epsilon}^1 t^{-d/p - \epsilon} \ \phi(t) \ dt, \ \epsilon \in (0, 1/2), \eqno(3.4)
$$

$$
K = \sup_{\psi} \sup_{ f \in G\psi, f \ne 0} \left[ \frac{||U_{\phi}[f]||G\psi_{\theta}}{||f||G\psi}   \right];
$$
it remains to prove $  K \ge 1. $

 Notice that

 $$
 K = \sup_{\psi} \sup_{ f \in G\psi, f \ne 0} \frac{[\sup_p |U_{\phi}[f]_p/(\theta(p) \psi(p))] }{\sup_p[ |f|_p/\psi(p)]}.
 $$
We conclude choosing $ \psi(p) = |f|_p  $

$$
K \ge \sup_{ f \in G\psi, f \ne 0} \left\{ \frac{|U_{\phi}[f]|_p}{\theta(p) \ |f|_p} \right\}. \eqno(3.5)
$$

\vspace{3mm}

We have after some calculations  using at the same example as in \cite{Xiao1}

$$
f_{\epsilon}(x) = I(|x| \ge 1) \cdot |x|^{-d/p - \epsilon}, \ |x|= \sqrt{(x,x)}: \eqno(3.6)
$$

$$
|f_{\epsilon}|_p^p = \frac{c(d)}{p \epsilon}, \ c(d) = \frac{d  \pi^{d/2}}{\Gamma(1 + d/2)};  \eqno(3.7)
$$

$$
U_{\phi}[f_{\epsilon}] = I(|x| > 1) \cdot |x|^{-d/p - \epsilon} \cdot \int_{1/|x|}^1 t^{ - d/p - \epsilon} \ \phi(t) \ dt;
$$

$$
|U_{\phi} [f_{\epsilon}]|_p^p =
\int_{|x| \ge 1} \left( |x|^{-d/p - \epsilon} \int_{1/|x|}^1 t^{-d/p - \epsilon} \ \phi(t) \ dt \right)^p \ dx \ge
$$

$$
|f_{\epsilon}|_p^p \cdot \left( \Delta^{-\epsilon} \int_{1/\Delta}^1 t^{-d/p} \ \phi(t) \ dt \right)^p; \eqno(3.8)
$$

$$
|U_{\phi} [f_{\epsilon}]|_p  \ge |f_{\epsilon}|_p \cdot\theta_{\epsilon}(p); \eqno(3.9)
$$

$$
K \ge \sup_{\epsilon \in (0, 1/2)} \sup_p \left[ \frac{\theta_{\epsilon}(p) \ |f_{\epsilon}|_p}{\theta(p) \ |f_{\epsilon}|_p}\right] \ge
\sup_p \overline{\lim}_{\epsilon \to 0+} \frac{\theta_{\epsilon}(p)}{\theta(p)}   = 1. \eqno(3.10)
$$

\hfill $\Box$

\bigskip

\section{Multidimensional case}

\vspace{3mm}

 We recall here the definition of the so-called anisotropic Lebesgue (Lebesgue-
Riesz) spaces, which arrear in the famous article belonging to Benedek A. and Panzone R.
\cite{Benedek1}. More detail information about this spaces see in the books of Besov O.V.,
Ilin V.P., Nikolskii S.M. \cite{Besov1}, chapter 16,17; Leoni G. \cite{Leoni1}, chapter 11.

  Let $ (X_j ,A_j , \mu_j, \ j = 1,2,\ldots,d) $ be measurable spaces with sigma-finite non - trivial measures
$ μ_j . $   It is clear that in this article $ X_j = R^{d_j}  $ and $ \mu_j  $ is ordinary Lebesgue measure. \par

 Let  $ p = \vec{p} = (p_1, p_2, ..., p_d) $ be $ d \ − \ $ dimensional vector such that $ 1 \le p_j \le \infty. $
Recall that the {\it anisotropic} Lebesgue space $ L(\vec{p}) $ consists on all the total measurable
real valued function $ f = f(x_1, x_2, . . . , x_d) = f(x) = f(\vec{x}), \ x_j \in X_j $
with finite norm $ |f|_{\vec{p} } \stackrel{def}{=} $

$$
\left( \int_{X_d} \mu_d(dx_d) \left(  \int_{X_{d-1}} \mu_{d-1}(dx_{d-1}) \ldots
\left( \int_{X_1} \mu_1(dx_1) |f(x_1,x_2, \ldots, x_d)|^{p_1}  \right)^{p_2/p_1} \right)^{p_3/p_2}  \ldots \right)^{1/p_d}.
$$

 Note that in general case $ |f|{p_1,p_2} \ne |f|{p_2,p_1}, $ but $ |f|_{p,p} = |f|_p. $
Observe also that if $ f(x_1, x_2) = g_1(x_1) · g_2(x_2), $ (condition of factorization), then
$ |f|_{p_1,p_2} = |g_1|_{p_1} · |g_2|_{p_2}, $ (formula of factorization).\par

\vspace{3mm}

  It is obvious that in the multidimensional  (anisotropic) case the inequality of Xiao look as follows:

  $$
  | U_{\phi} [f]|_{\vec{p}} \le \prod_{k=1}^l \theta_{d_k}(p_k) \cdot |f|_{\vec{p}}. \eqno(4.1)
  $$
The last estimate is exact, for instance, for factorable function:

$$
f(\vec{x}) = \prod_{k=1}^k f_k( \vec{x_k}).
$$

\vspace{3mm}

{\bf Anisotropic Grand Lebesgue-Riesz spaces.}

\vspace{3mm}
 Let $ Q $ be convex (bounded or not) subset of the set $ \otimes_{j=1}^l [1,\infty]. $
 Let $ \psi = \psi(\vec{p}) $ be continuous in an interior $ Q^0 $ of the set $ Q $
strictly  positive  function such that

$$
\inf_{\vec{p} \in Q^0}  \psi(\vec{p}) > 0; \ \inf_{\vec{p} \notin Q^0}  \psi(\vec{p}) = \infty.
$$

 We denote the set all of such a functions as $ \Psi(Q). $ \par
The  Anisotropic Grand Lebesgue Spaces $ AGLS = AGLS(\psi) $ space consists on all the measurable functions

$$
f:  \otimes_{j=1}^l X_j \to R
$$
with finite (mixed) norms

$$
||f||AG\psi = \sup_{\vec{p} \in Q^0} \left[ \frac{|f|_{\vec{p}}}{\psi(\vec{p} )} \right].\eqno(4.2)
$$

 An application  into the theory of multiple Fourier transform of these spaces see in articles \cite{Benedek1} and
 \cite{Ostrovsky100}, where are considered some problems
of boundedness of singular operators  in (weight) Grand Lebesgue Spaces and in Anisotropic Grand Lebesgue Spaces.\par

\vspace{3mm}

 Let $ \psi \in \Psi(Q); $ we denote

 $$
 \psi_{\theta} (\vec{p}) = \prod_{k=1}^l \theta_{d_k}(p_k) \cdot \psi(\vec{p}). \eqno(4.3)
 $$

\vspace{3mm}

 We deduce analogously to the theorem 3.1:\\

\vspace{3mm}

{\bf  Proposition 4.1. }\\

\vspace{3mm}

$$
||U_{\phi} [f]||G\psi_{\theta} \le 1 \cdot ||f||G\psi, \eqno(4.4)
$$
\vspace{3mm}
{\it where the constant "1" is the best possible. } \par

\hfill $\Box$

\bigskip

\section{Concluding remarks}

\vspace{3mm}

{\bf 1.  Analogously may be investigated the "conjugate" operator  of a view}

 $$
 V_{\phi}[f](x) = \int_0^1 f(x/t) \ t^{-d} \ \phi(t) \ dt.
 $$

 J.Xiao in \cite{Xiao1}  proved that

$$
||V_{\phi}||(L(p) \to L(p)) = \zeta(p) \stackrel{def}{=} \int_0^1 t^{-d(1 - 1/p)} \ \phi(t) \ dt.
$$

 Define for any function $ \psi = \psi(p) $ a new function

 $$
 \psi^{(\zeta)}(p) = \psi(p) \ \zeta(p).
 $$
We conclude alike to the assertion of theorem 3.1:

$$
||V_{\phi}[f]|| G\psi^{(\zeta)} \le 1 \cdot ||f||G\psi,
$$
with exact value of the coefficient "1". \par

\vspace{4mm}

{\bf 2. Examples.}\par

\vspace{3mm}

{\bf A.} Let

$$
\phi(t) =  t^{\alpha - 1} \ (1-t)^{\beta - 1}, \ \alpha,\beta = \const > 0;
$$
then

$$
\theta(p)  = B(\alpha - d/p, \beta) = \frac{\Gamma (\alpha - d/p) \ \Gamma(\beta)}{ \Gamma(\alpha + \beta - d/p)}, \ p > d/\alpha
$$
and $ \theta(p) = \infty  $ otherwise. Here as usually $ B(\cdot, \cdot), \ \Gamma(\cdot) $ denote correspondingly Beta and Gamma
functions.\par
 Note that as $ p \to d/\alpha +0, \ p > d/\alpha   \Rightarrow \theta(p) \sim  p/(\alpha p - d). $ \\

\vspace{3mm}

{\bf B.}  Let now

$$
\phi(t) = |\log t|^{\gamma} \ L(|\log t|), \ \gamma = \const > -1,
$$
$ L  = L(z), \ z \in (0, \infty) $ is positive continuous slowly varying as $ z \to \infty $ function.
  We have as $ p \to d+0, \ p > d $

$$
\theta(p) = \int_0^{1} t^{-d/p} \ |\log t|^{\gamma} \ L(|\log t|) \ dt = \int_0^{\infty} e^{-y(1 - d/p)} \ y^{\gamma} \ L(y) dy =
$$

$$
 \left[\frac{p}{p-d} \right]^{\gamma + 1} \int_0^{\infty} e^{-z} z^{\gamma} L \left( \frac{pz}{p-d}   \right)\ dz \sim
$$

$$
\Gamma(\gamma + 1) \left[\frac{p}{p-d} \right]^{\gamma + 1} L \left( \frac{p}{p-d}   \right).
$$

 Evidently, $ \theta(p) = \infty, \ p \le d. $\\

\vspace{3mm}

{\bf  C.  Multivariate Hardy transform.} \par

\vspace{3mm}

In 1976 Faris \cite{Faris1} first gave a definition of Hardy operator in $ d - $ dimensional case. In
1995, Christ and Grafakos \cite{Christ1}  gave its equivalent version of $  d - $ dimensional Hardy operator
as follows

$$
H_d(f)(x) = \frac{1}{\Omega(d) |x|^d} \int_{y: |y| \le |x|} f(y) dy, \ \Omega(d) = \frac{\pi^{d/2}}{\Gamma(1 + d/2)}, \ x \ne 0.
$$
and proved that

$$
|| H_d ||(L(p) \to L(p)) = \frac{p}{p-1}, \ p > 1,
$$
as in the one - dimensional case.\par
 Define the following transform for every $ \psi - $ function:

 $$
 \psi_1(p) = \psi(p) \ \frac{p}{p-1}, \ p > 1;
 $$
then

$$
||H_d [f]||G\psi_1 \le 1 \cdot ||f||G\psi,
$$
where the constant "1" is non-improvable.\par

\end{document}